\documentclass[12pt,a4paper,oneside]{amsart}

\usepackage{amscd}
\input xy
\xyoption{all}

\newtheorem{definition}{Definition}
\newtheorem{thm}{Theorem}
\newtheorem{lem}{Lemma}
\newtheorem{prop}{Proposition}
\newtheorem{cor}{Corollary}

\newcommand{\To}{\longrightarrow}
\newcommand{\Complex}{\mathbb C}

\DeclareMathOperator{\End}{End}
\DeclareMathOperator{\Id}{Id}
\DeclareMathOperator{\Hom}{Hom}
\DeclareMathOperator{\Sym}{Sym}

\DeclareMathOperator{\tr}{tr}
\DeclareMathOperator{\sign}{sign}
\DeclareMathOperator{\Arrows}{Arrows}
\DeclareMathOperator{\Mat}{Mat}
\DeclareMathOperator*{\Ast}{\times}
\numberwithin{equation}{section}
\numberwithin{prop}{section}
\numberwithin{cor}{section}
\numberwithin{lem}{section}

\begin{document}
\title[Algebras generated by elements with given spectrum\dots]{Algebras generated by elements with given spectrum and scalar sum and
Kleinian singularities}
\author{Anton Mellit}

\email{mellit@imath.kiev.ua}
\address{Institute of Mathematics, National Academy of Sciences of
Ukraine, 3 Tereshchenkivska Street, Kyiv 4, 01601, Ukraine}

\begin{abstract}
We consider algebras $e_i \Pi^\lambda(Q) e_i$ where $\Pi^\lambda(Q)$ is the deformed 
preprojective algebra of weight $\lambda$ and $i$ is some vertex of $Q$ in the case
when $Q$ is an extended Dynkin diagram and $\lambda$ lies on the hyperplane
orthogonal to the minimal positive imaginary root $\delta$. We prove that the center
of $e_i \Pi^\lambda(Q) e_i$ is isomorphic to $\mathcal{O}^\lambda(Q)$ which
is a deformation of coordinate ring of Kleinian singularity which corresponds to $Q$.
Also we find the minimal $k$ for which the standard identity of degree $k$
holds in $e_i \Pi^\lambda(Q) e_i$. We prove that algebras $A_{P_1,\dots,P_n;\mu} = 
\mathbb{C}\langle x_1, \dots, x_n | P_i(x_i)=0,
    \sum_{i=1}^n x_i = \mu e\rangle$ are the special case of algebras 
$e_c \Pi^\lambda(Q) e_c$ for star-like quivers $Q$ with origin $c$.
\end{abstract}

\maketitle
\section*{Introduction}

Consider the problem of description of n-tuples of hermitian operators $\{A_i\}$
in a Hilbert space satisfying given restrictions on spectra $\sigma(A_i)\subset M_i$ with
$M_i\subset\mathbb{R}$ finite
and relation $\sum_{i=1}^n A_i = \mu I$, with $I$ - the identity and $\mu\in \mathbb{R}$.
Study of such n-tuples is equivalent to study of *-representations of certain *-algebra.
Forgetting the *-structure we arrive to the following class of algebras.

\begin{definition}
Let $P_1$, \dots, $P_n$ be complex polynomials in one variable and $\mu\in \mathbb{C}$.
We put inessential restriction $P_i(0)=0$. Define algebra
\[
    A_{P_1,\dots,P_n;\mu} = \mathbb{C}\langle x_1, \dots, x_n | P_i(x_i)=0\, (i=1,\dots,n),
    \sum_{i=1}^n x_i = \mu e\rangle.
\]
\end{definition}

In joint work of the author with Yu. Samoilenko and M. Vlasenko (see \cite{ours})
we studied some properties of such algebras:
we computed growth of these algebras and proved existence of polynomial identities in
certain cases (in fact, finiteness over center was proved).

These algebras are closely related to deformed preprojective algebras of W.~Crawley-Boevey
and M.P.~Holland (\cite{CrB}). We briefly recall their definition. Let $Q$ be a quiver with
vertex set $I$. Write $\bar{Q}$ for the double quiver of $Q$, i.e. quiver obtained by adding 
a reverse arrow $a^*:j\To i$ for every arrow $a:i\To j$, and write $\Complex \bar{Q}$ for its
path algebra, which has basis the paths in $\bar{Q}$, including a trivial path $e_i$  for each 
vertex $i$. If $\lambda=(\lambda_i) \in \Complex^I$, then the deformed preprojective algebra of
weight $\lambda$ is
\[
\Pi^\lambda(Q) = \Complex \bar{Q} / (\sum_{a\in \Arrows(Q)}[a,a^*]-\lambda),
\]
where $\Arrows(Q)$ denotes the set of arrows of $Q$, and $\lambda$ is 
identified with the element $\sum_{i\in I} \lambda_i e_i$.

Let $A = A_{P_1,\dots,P_n;\mu}$. Consider quiver $Q(A)$ with vertices
\[
I=\{(i, j) | i=1,\dots,n ,\, j=1,\dots,\deg P_i-1\} \cup \{c\}
\]
and arrows
\[
\{a_{ij} : (i, j) \To (i,j-1) | i=1,\dots,n ,\, j=1,\dots,\deg P_i-1\},
\]
where $(i, 0)$ is identified with $c$ for $i=1,\dots,n$.
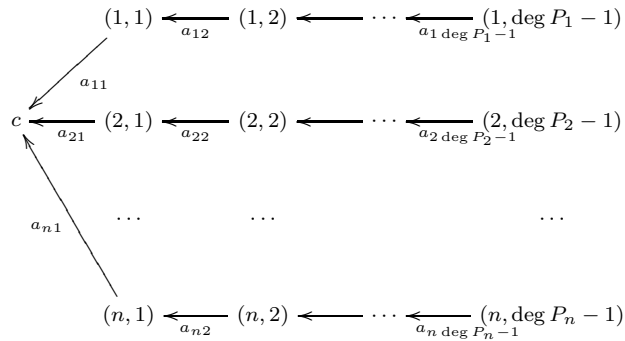
\[
\begin{tiny}
	\xymatrix{
	&(1, 1)\ar@{->}[ld]^{a_{11}} & (1, 2)\ar@{->}[l]^{a_{12}} & \cdots \ar@{->}[l] & (1, \deg P_1-1) \ar@{->}[l]^{a_{1\deg P_1-1}}\\
c	&(2, 1)\ar@{->}[l]^{a_{21}} & (2, 2)\ar@{->}[l]^{a_{22}} & \cdots \ar@{->}[l] & (2, \deg P_2-1) \ar@{->}[l]^{a_{2\deg P_2-1}}\\
	&\cdots & \cdots&&\cdots\\
	&(n, 1)\ar@{->}[luu]^{a_{n1}} & (n, 2)\ar@{->}[l]^{a_{n2}} & \cdots \ar@{->}[l] & (n, \deg P_n-1) \ar@{->}[l]^{a_{n\deg P_n-1}}\\
	}
\end{tiny}	
\]

Note that the graph $Q$ coincides with the graph of algebra $A$, considered in
\cite{ours}. Here is an example of quiver $Q$ for the case $n=3$, $\deg P_1 = 2$, $\deg P_2 = 3$,
$\deg P_3 = 2$:
\[
\begin{CD}
e_{11} @>{a_{11}}>> e_c @<{a_{21}}<< e_{21} @<{a_{22}}<< e_{22}\\
@. @AA{a_{31}}A\\
@. e_{31}
\end{CD}
\]

The first result establishes connection between algebras $A_{P_1,\dots,P_n;\mu}$ and
deformed preprojective algebras.
\begin{thm}\label{thm1}
Algebra $A=A_{P_1,\dots,P_n;\mu}$ is isomorphic to $e_c\Pi^{\lambda}(Q)e_c$ under 
the isomorphism sending $x_i$ to $a_{i1} a_{i1}^*$ for $Q=Q(A)$ and
\[
\lambda = \sum_{i=1}^n\sum_{j=1}^{\deg P_i-1} (\alpha_{ij-1}-\alpha_{ij})e_{ij}+\mu e_c,
\]
where $\alpha_{i0}, \alpha_{i1}, \dots, \alpha_{i \deg P_i-1}$ are all
roots of the polynomial $P_i$ taken with multiplicities in any order with
 $\alpha_{i0}=0$.
\end{thm}

Consider the case when graph $Q$ is an extended Dynkin diagram of type $\widetilde{A_n}$, 
$\widetilde{D_n}$ or $\widetilde{E_n}$. The following picture shows all such graphs along with
coordinates of the so called minimal imaginary root $\delta\in \Complex^I$. Boxed vertex is
the extending vertex.
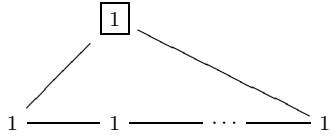
\begin{description}
\item[$\widetilde{A_n}$]	\[\begin{tiny}
		\xymatrix{
		&\boxed{1} \ar@{-}[ld]&&\\
		1\ar@{-}[r] & 1\ar@{-}[r] & \cdots\ar@{-}[r] & 1\ar@{-}[ull]\\
		}
	\end{tiny}\]	
\item[$\widetilde{D_n}$]	
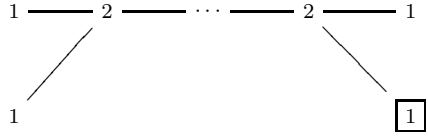
\[\begin{tiny}
		\xymatrix{
		1\ar@{-}[r] & 2\ar@{-}[r] & \cdots\ar@{-}[r] & 2\ar@{-}[r] & 1\\
		1\ar@{-}[ru]&&&& \boxed{1}\ar@{-}[lu]\\
		}
	\end{tiny}\]
\item[$\widetilde{E_6}$]	
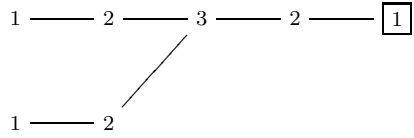
\[\begin{tiny}
		\xymatrix{
		1\ar@{-}[r] & 2\ar@{-}[r] & 3 & 2\ar@{-}[l] & \boxed{1}\ar@{-}[l]\\
		1\ar@{-}[r] & 2\ar@{-}[ru]&&&\\
		}
	\end{tiny}\]
\item[$\widetilde{E_7}$]	
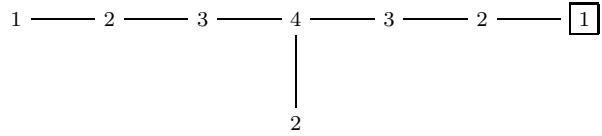
\[\begin{tiny}
		\xymatrix{
		1\ar@{-}[r] & 2\ar@{-}[r] & 3\ar@{-}[r] & 4 & 3\ar@{-}[l] & 2\ar@{-}[l] & \boxed{1}\ar@{-}[l]\\
		& && 2\ar@{-}[u]&&&\\
		}
	\end{tiny}\]
\item[$\widetilde{E_8}$]	
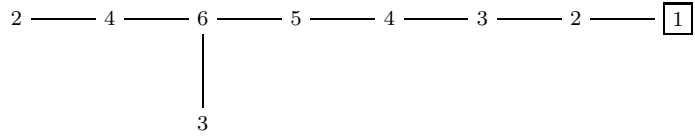
\[\begin{tiny}
		\xymatrix{
		2\ar@{-}[r] & 4\ar@{-}[r] & 6 & 5\ar@{-}[l] & 4\ar@{-}[l] & 3\ar@{-}[l] & 2\ar@{-}[l] & \boxed{1}\ar@{-}[l]\\
		& & 3\ar@{-}[u]&&&&&\\
		}
	\end{tiny}\]
\end{description}

In \cite{CrB} they proved that $\Pi^{\lambda}(Q)$ is a PI-algebra 
(on PI algebras see \cite{Rowen1990})
if and only if $\delta \cdot \lambda=0$.
Also they studied algebra $\mathcal{O}^{\lambda}(Q)$ which is $e_0 \Pi^{\lambda}(Q) e_0$ where
$0$-th vertex is the extending vertex of $Q$, and proved that this algebra is commutative if and
only if $\delta \cdot \lambda=0$.
For $\lambda=0$ the algebra $\mathcal{O}^0(Q)$ coincides with the coordinate ring
of the corresponding Kleinian singularity.

In this paper we consider algebras $e_i \Pi^{\lambda} e_i$ for arbitrary $i\in I$. For the case
$\delta \cdot \lambda=0$ we study center of such algebra and find minimal number $k$ for which 
it possesses standard identity of degree $k$, i.e.
\[
\sum_{\pi\in \mathcal{S}_k} \sign(\pi) \prod_{i=1}^k x_{\pi(i)} = 0.
\]
We denote by $\mathcal{S}_k$ the group of permutations of $k$ elements.

We obtain the following theorems:

\begin{thm}\label{thmCent}
If $Q$ is an extended Dynkin diagram $\widetilde{A_n}$, 
$\widetilde{D_n}$ or $\widetilde{E_n}$, $\delta \cdot \lambda=0$ and $i\in I$ is some vertex of $Q$ then
center of $e_i \Pi^{\lambda}(Q) e_i$ is isomorphic to $\mathcal{O}^{\lambda}(Q)= e_0 \Pi^{\lambda}(Q) e_0$ where
$0$-th vertex is the extending vertex of $Q$.
\end{thm}

\begin{thm}\label{thmPI}
If $Q$ is an extended Dynkin diagram $\widetilde{A_n}$, 
$\widetilde{D_n}$ or $\widetilde{E_n}$, $\delta \cdot \lambda=0$ and $i\in I$ is some vertex of $Q$ then
$e_i \Pi^{\lambda} e_i$ possesses standard identity of degree $2 \delta_i$ and it is the minimal number
with such property.
\end{thm}

\section{Representations of groups}
Suppose $V$ is a two dimensional complex vector space with simplectic form $\omega$.
Let $G$ be a finite subgroup of $SL(V)$. Suppose irreducible representations of $G$ are
precisely $\{V_i\}_{i\in I}$ where $I=\{0,1,2,\dots, n\}$ with $V_0$~--- the trivial one. Suppose 
\[
	V \otimes V_i = \bigoplus_{j=1}^n m_{ij} V_j.
\]
Then the McKay graph of $G$ is defined to be a graph with vertex set $I$ and number
of edges between $i$ and $j$ is $m_{ij}$ (we will always have $m_{ij}=m_{ji}$).
According to J. McKay (\cite{McKay}) McKay graphs of finite subgroups of $SL(V)$ are
extended Dynkin diagrams: $\widetilde{A_n}$ for cyclic groups, $\widetilde{D_n}$ for dihedral groups
and $\widetilde{E_6}$, $\widetilde{E_7}$, $\widetilde{E_8}$ for binary tetrahedral, octahedral and 
icosahedral groups. Dimensions of irreducible representations are $\dim V_i = \delta_i$.
Let us fix some orientation of McKay graph of $G$ thus obtaining a quiver $Q$.

Let $M$ be some $\Complex G$-module. Consider vector space $F_0(M) = T(V^*)\otimes M$ where 
\[
T(V^*) = \bigoplus_{i=0}^{\infty} V^{*\otimes i} \; \text{--- the tensor algebra of $V^*$}.
\]
Equip $F_0(M)$ with the componentwise action of $G$. Then we can consider $F(M)=F_0(M)^G$~--- 
the subspace of $G$-invariant vectors. Note that if $M$ is algebra and multiplication respects
action of $G$ then both $F_0(M)$ and $F(M)$ become graded algebras with grading 
$F_0(M)_i = V^{*\otimes i} \otimes M$ and $F(M)_i = (F_0(M)_i)^G$. Consider algebra 
$F(\End_\Complex(V_{\Sigma}))$, where $V_{\Sigma}$ is the direct sum of all irreducible $\Complex G$-modules and 
$G$ acts on $\End_\Complex(V_{\Sigma})$ by conjugation.
Clearly $F(\End_\Complex(V_{\Sigma}))_0$ is the same as $(\End_\Complex(V_\Sigma))^G$, which 
in its turn can be identified with $\Complex^I$. Our aim is to build a graded algebra isomorphism $\varphi$
from $\Complex \bar{Q}$ to $F(\End_\Complex(V_{\Sigma}))$ such that 
\begin{equation}\tag{*}
\begin{split}
\varphi \text{ is identity on } \Complex^I \\
\varphi(\sum_{a\in \Arrows(Q)} [a, a^*]) = \delta \omega.
\end{split}
\end{equation}
We accomplish this in two steps:
\begin{lem}\label{lemTens}
Natural homomorphism 
\[
F(\End_\Complex(V_{\Sigma}))_i \otimes_{F(\End_\Complex(V_{\Sigma}))_0} F(\End_\Complex(V_{\Sigma}))_j \To
F(\End_\Complex(V_{\Sigma}))_{i+j}
\]
is an isomorphism
\end{lem}
\begin{proof}
We make some identifications
\[
\begin{split}
	F(\End_\Complex(V_{\Sigma}))_i = (V^{*\otimes i} \otimes \End_\Complex(V_{\Sigma}))^G \cong
	\Hom_G(V_\Sigma, V^{*\otimes i} \otimes V_\Sigma) \\
	\cong \bigoplus_{i\in I} \Hom_\Complex(\Complex, \Complex^{a_i}),
\end{split}
\]	
where 
\[
V^{*\otimes i} \otimes V_\Sigma \cong \bigoplus_{i\in I} V_i^{\oplus a_i}.
\]
\[
\begin{split}
	F(\End_\Complex(V_{\Sigma}))_j = (V^{*\otimes j} \otimes \End_\Complex(V_{\Sigma}))^G \cong
	\Hom_G(V^{\otimes j}\otimes V_\Sigma, V_\Sigma) \\
	\cong \bigoplus_{i\in I} \Hom_\Complex(\Complex^{b_i}, \Complex),
\end{split}
\]
where
\[
V^{\otimes j}\otimes V_\Sigma \cong \bigoplus_{i\in I} V_i^{\oplus b_i}.
\]
\[
\begin{split}
	F(\End_\Complex(V_{\Sigma}))_{i+j} = (V^{*\otimes (i+j)} \otimes \End_\Complex(V_{\Sigma}))^G \\
	\cong \Hom_G(V^{\otimes j}\otimes V_\Sigma, V^{*\otimes i} \otimes V_\Sigma)
	\cong \bigoplus_{i\in I} \Hom_\Complex(\Complex^{b_i}, \Complex^{a_i}).
\end{split}
\]	
Recall that 
\[
F(\End_\Complex(V_{\Sigma}))_0 \cong \bigoplus_{i\in I} \Complex.
\]
Now the statement is clear.
\end{proof}

This lemma implies that the natural homomorphism from tensor algebra of 
$F(\End_\Complex(V_{\Sigma}))_1$ over $F(\End_\Complex(V_{\Sigma}))_0$ to 
$F(\End_\Complex(V_{\Sigma}))$ is an isomorphism. Graded algebra $\Complex \bar{Q}$ 
possesses the same property, so it is clear that for establishing isomorphism of graded algebras
which is identity on $\Complex^I$ it is necessary and sufficient to establish 
isomorphism of subbimodules in degree $1$. Decompose $F(\End_\Complex(V_{\Sigma}))_1$ by primitive
idempotents of $\Complex^I$:
\[
F(\End_\Complex(V_{\Sigma}))_1 = (V^* \otimes \End_\Complex(V_{\Sigma}))^G \cong
\bigoplus_{i,j\in I} \Hom_G(V\otimes V_i, V_j).
\]
Clearly $\Hom_G(V\otimes V_i, V_j)$ vanish if there are no arrows from $i$ to $j$ in $\bar{Q}$ and
is one dimensional if there is arrow from $i$ to $j$ in $\bar{Q}$. Subbimodule of $\Complex^I$ in degree
$1$ has similar decomposition. So any assignment $a\To \varphi(a) \in \Hom_G(V\otimes V_i, V_j)$, 
$\varphi(a) \neq 0$ for $a\in \Arrows(\bar{Q})$ induces some isomorphism of graded algebras 
$\varphi: \Complex \bar{Q} \To F(\End_\Complex(V_{\Sigma}))$. 
\begin{prop}
For every arrow $a:i\To j$ of $Q$ choose any nonzero representative $\varphi(a)\in \Hom_G(V\otimes V_i, V_j)$.
It is possible to choose $\varphi(a^*)\in \Hom_G(V_j, V^* \otimes V_i)$ such that 
\[
\tr((\iota \otimes \Id_{V_i}) \varphi(a^*) \varphi(a) ) = \dim V_i \dim V_j,
\]
where $\iota: V^* \To V$ is such that 
\[
f(x) = \omega(\iota(f), x) \; \text{for $f\in V^*$ and $x\in V$}.
\]
This induces isomorphism of algebras which satisfies property (*).
\end{prop}
\begin{proof}
As for the possibility of choosing such a $\varphi(a^*)$, in decomposition 
of $V\otimes V_i$ into direct sum of indecomposable $\Complex G$-modules $V_j$ occurs
exactly once, so if we choose any nonzero $\varphi(a^*)\in \Hom_G(V\otimes V_i, V_j)$, we 
obtain that 
\[
(\iota \otimes \Id_{V_i}) \varphi(a^*) \varphi(a)
\]
is a projection on $V_j$ in $V\otimes V_i$ multiplied by some complex constant, so its trace is 
nonzero and by multiplication by some factor it is possible to make the trace accepting any complex
value.
It is only needed check that 
\[
\sum_{a\in \Arrows(Q)} [\varphi(a), \varphi(a^*)] = \delta \omega.
\]
Choose some vertex $i$ and multiply both sides by $e_i$:
\begin{equation} \label{eq10}
\sum_{\substack{j\in I,\, a: j\To i,\\ a\in \Arrows(Q)}} \varphi(a) \varphi(a^*) - 
\sum_{\substack{j\in I,\, a: i\To j,\\ a\in \Arrows(Q)}} \varphi(a^*) \varphi(a) = \delta_i \omega e_i.
\end{equation}
Both sides belong to $\Hom_G(V \otimes V \otimes V_i, V_i)$, which can be 
identified with $\Hom_G(V \otimes V_i, V^* \otimes V_i)$ by 'lifting' first element of tensor 
product. Apply $\iota \otimes \Id_{V_i}$ to both sides. Since $(\omega(x))(y) = \omega(y, x)$ and
$(\omega(x))(y) = \omega(\iota(\omega(x)), y)$ we have $\iota(\omega(x)) = -x$ and
\[
(\iota \otimes \Id_{V_i}) \delta_i \omega e_i = -\delta_i \Id_{V \otimes V_i}.
\]
Recall that each $(\iota \otimes \Id_{V_i})\varphi(a) \varphi(a^*)$ and
$(\iota \otimes \Id_{V_i})\varphi(a^*) \varphi(a)$ which occurs in (\ref{eq10}) is a projection on a summand $V_j$ 
multiplied by some complex number where $j$ is another endpoint of $a$ different from $i$. 
Denote this projection by $p_j$. Then
\[
(\iota \otimes \Id_{V_i})\varphi(a) \varphi(a^*) = 
	\frac{\tr((\iota \otimes \Id_{V_i})\varphi(a) \varphi(a^*))}{\dim V_j} p_j \;\; \text{and}
\]
\[
-(\iota \otimes \Id_{V_i})\varphi(a^*) \varphi(a) = 
	\frac{-\tr((\iota \otimes \Id_{V_i})\varphi(a^*) \varphi(a))}{\dim V_j} p_j.
\]
By definition 
\[
\tr((\iota \otimes \Id_{V_i})\varphi(a^*) \varphi(a)) = \dim V_i  \dim V_j.
\]
There is an identity
\[
\tr((\iota \otimes \Id_{V_i})x y) = 
-\tr((\iota \otimes \Id_{V_j})y x),
\]
which holds for every $x\in \Hom_\Complex(V \otimes V_j, V_i)$ and $y\in \Hom_\Complex(V \otimes V_i, V_j)$.
It is enough to check this identity for $x = f_1 \otimes x_0$ and $y = f_2 \otimes y_0$ 
where $f_1, f_2 \in V^*$, $x_0 \in \Hom_\Complex(V_j, V_i)$ and $y_0 \in \Hom_\Complex(V_i, V_j)$:
\begin{multline*}
\tr((\iota \otimes \Id_{V_i})x y) = \tr(\iota(f_1)f_2 \otimes x_0 y_0) = f_2(\iota(f_1)) \tr(x_0 y_0) \\= 
\omega(\iota(f_2), \iota(f_1)) \tr(x_0 y_0) = -\omega(\iota(f_1), \iota(f_2)) \tr(y_0 x_0) \\=
-f_1(\iota(f_2)) \tr(y_0 x_0) = -\tr(\iota(f_2)f_1 \otimes y_0 x_0) = 
-\tr((\iota \otimes \Id_{V_j})y x).
\end{multline*}
Apply this identity:
\[
\tr((\iota \otimes \Id_{V_i})\varphi(a) \varphi(a^*)) = -\tr((\iota \otimes \Id_{V_j})\varphi(a^*) \varphi(a)) = -\dim V_i  \dim V_j.
\]
It follows that $\iota \otimes \Id_{V_i}$ applied to lefthand side of (\ref{eq10}) equals to 
\[
-\dim V_i \sum_{\substack{j\in I,\, a: i\To j,\\ a\in \Arrows(\bar{Q})}} p_j = -\dim V_i \Id_{V\otimes V_i}
\]
and recalling that $\delta_i = \dim V_i$ we are done.
\end{proof}

The next corollary summarizes what have been done in this section.
\begin{cor}
Algebra $\Pi^\lambda(Q)$ is isomorphic to algebra
\[
(T(V^*)\otimes \End_\Complex(V_\Sigma))^G / (\delta \omega - \lambda).
\]
Moreover, this is isomorphism of filtered algebras with filtrations induced from
gradings of $\Complex \bar{Q}$ and $T(V^*)$.
\end{cor}

\section{Case of $\lambda=0$}
In this section we are going to prove theorems \ref{thmCent} and \ref{thmPI} for the 
case $\lambda=0$. Key is the following lemma:
\begin{lem}
Algebra $\Pi^0(Q)$ is isomorphic to algebra of polynomial $G$-equivariant maps 
from $V$ to $\End_\Complex(V_\Sigma)$, i.e. the algebra
\[
(\Sym(V^*)\otimes \End_\Complex(V_\Sigma))^G,
\]
where $\Sym(V^*)$ is the algebra of symmetric tensors of $V^*$. Moreover, this is isomorphism
of graded algebras. 
\end{lem}
\begin{proof}
We already know that $\Pi^0(Q)$ is isomorphic to 
\[
(T(V^*)\otimes \End_\Complex(V_\Sigma))^G / (\delta \omega) = (T(V^*)\otimes \End_\Complex(V_\Sigma))^G / \omega.
\]
Since
\[
\Sym(V^*)\otimes \End_\Complex(V_\Sigma) = (T(V^*) / w)\otimes \End_\Complex(V_\Sigma) =
(T(V^*) \otimes \End_\Complex(V_\Sigma)) / w
\]
it is enough to prove that the idempotent 
\[
\varepsilon = \frac{1}{|G|} \sum_{g\in G} g
\]
maps ideal generated by $\omega$ in $T(V^*) \otimes \End_\Complex(V_\Sigma)$ to the ideal
generated by $\omega$ in $(T(V^*) \otimes \End_\Complex(V_\Sigma))^G$. To prove this take
some $f \in V^{*\otimes i} \otimes \End_\Complex(V_\Sigma)$, 
$g \in V^{*\otimes j} \otimes \End_\Complex(V_\Sigma)$ and consider $\varepsilon (f \omega g)$.
Note that $f\omega g$ is antisymmetric with respect to $i+1$-th and $i+2$-th argument. It follows
that $\varepsilon(f \omega g)$ is antisymmetric with respect to $i+1$-th and $i+2$-th argument
as well. Since 
\[
\varepsilon(f \omega g) \in (V^{*\otimes (i+j+2)}\otimes \End_\Complex(V_\Sigma))^G
\]
and we know from lemma \ref{lemTens} that 
\[
\begin{split}
(V^{*\otimes (i+j+2)}\otimes \End_\Complex(V_\Sigma))^G \\= 
(V^{*\otimes i}\otimes \End_\Complex(V_\Sigma))^G \otimes_{\Complex^G}
(V^{*\otimes 2}\otimes \End_\Complex(V_\Sigma))^G \otimes_{\Complex^G}
(V^{*\otimes j}\otimes \End_\Complex(V_\Sigma))^G
\end{split}
\]
we can decompose
\[
\varepsilon(f \omega g) = \sum_{k=1}^K f_k \omega_k g_k
\]
with $f_k \in (V^{*\otimes i}\otimes \End_\Complex(V_\Sigma))^G$, 
$g_k \in (V^{*\otimes j}\otimes \End_\Complex(V_\Sigma))^G$ and
$\omega_k \in (V^{*\otimes 2}\otimes \End_\Complex(V_\Sigma))^G$.
Denote by $\tau$ operator acting on elements of $(V^{*\otimes (i+j+2)}\otimes \End_\Complex(V_\Sigma))^G$
by interchanging $i+1$-th and $i+2$-th arguments. Then
\[
\tau \varepsilon(f\omega g) = \sum_{k=1}^K f_k \omega'_k g_k
\]
with $\omega'_k$ is obtained from $\omega_k$ by interchanging first two arguments.
Hence
\[
\varepsilon(f \omega g) = \frac12 (\varepsilon(f \omega g) - \tau \varepsilon(f \omega g)) = 
\frac12 \sum_{k=1}^K f_k (\omega_k-\omega'_k) g_k.
\]
Since $\omega_k-\omega'_k\in \Hom_G(V\otimes V, \End_\Complex(V_\Sigma))$ is antisymmetric and $V$ is
two dimensional, it can be represented
as $\omega x_k$ with $x_k\in \End_\Complex(V_\Sigma)^G$. Thus
\[
\varepsilon(f \omega g) = \frac12 \sum_{k=1}^K f_k \omega x_k g_k
\]
with $f_k$, $x_k$ and $g_k$ from $(T(V^*) \otimes \End_\Complex(V_\Sigma))^G$. This completes the proof.
\end{proof}

The next propositions follow immediately.
\begin{prop}
Algebra $e_i\Pi^0(Q) e_i$ is isomorphic to the algebra of polynomial $G$-equivariant maps 
from $V$ to $\End_\Complex(V_i)$ for any $i\in I$. In particular, 
$\mathcal{O}^0(Q) = e_0 \Pi^0(Q) e_0$ is isomorphic to the algebra of invariants of $G$ on $V$.
\end{prop}
\begin{prop}
Algebra $e_i\Pi^0(Q) e_i$ possesses standard identity of degree $2 \delta_i$ for any $i\in I$.
\end{prop}
\begin{prop}\label{prop0Cent}
There is a graded inclusion from $e_0\Pi^0(Q) e_0$ to the center of $\Pi^0(Q)$ and 
graded inclusions from $e_0\Pi^0(Q) e_0$ to center of $e_i\Pi^0(Q) e_i$ for $i\in I$
induced by inclusions $\Complex \subset \End_\Complex(V_\Sigma)$ and $\Complex \subset \End_\Complex(V_i)$
correspondingly.
\end{prop}

For any $i\in I$ and $x\in V$ denote by $\mu_i(x)$ the subset of $\End_\Complex(V_i)$ defined by
\[
\mu_i(x) = \{f(x) | f \;\text{is a polynomial $G$-equivariant map from $V$ to $\End_\Complex(V_i)$}\}.
\]
For the rest we need the following statement:
\begin{lem}
The set of $x\in V$ for which $\mu_i(x) = \End_\Complex(V_i)$ is algebraically dense 
for any $i\in I$.
\end{lem}
\begin{proof}
Suppose $f: V\To \Complex$ is a non-constant $G$-invariant polynomial function. Then its differential $df$
is a polynomial $G$-equivariant map from $V$ to $V^*$. Denote by $U$ the set of $x\in V$ for which
$(df(x))(x) \neq 0$. Clearly $U$ is open and $U$ is not empty since $(df(x))(x)=0$ implies $f$ is a constant.
Denote by $U'$ the subset of $U$ of all $x$ such that $f(x)\neq 0$. Since $U'$ is open and not empty it is dense.
We will prove that every $x$ from $U'$ satisfies the required condition. So let $f(x)\neq 0$ and let 
$(df(x))(x) \neq 0$. Then $\iota df(x) \in V$ ($\iota: V^* \To V$ is such that 
$\omega(\iota(y_1), y_2) = y_1(y_2)$ for every $y_1\in V^*$ and $y_2\in V$) is not a multiple of $x$ because
if $\iota df(x) = C x$, $C \in \Complex$ then 
\[
(df(x))(x) = \omega(\iota df(x), x)=\omega(C x, x) = 0.
\]
It follows that $f(x) x$ and $\iota df(x)$ span $V$. Since $g_1 : V \To V$ defined by $g_1(y) = f(y) y$ and
$g_2 : V \To V$ defined by $g_2(y) = df(y)$ are polynomial and $G$-equivariant we have that every element
of $V$ is value in $x$ of some polynomial $G$-equivariant map from $V$ to $V$. It follows that for every
$k$ every element of $V^{\otimes k}$ is value in $x$ of some polynomial $G$-equivariant map from $V$ to 
$V^{\otimes k}$. Since every finite dimensional $\Complex G$-module is submodule of $V^{\otimes k}$ for some $k$, 
the statement holds for every finite dimensional $\Complex G$-module, in particular for $\End_\Complex(V_i)$.
\end{proof}

This lemma implies that there is no $k<2 \delta_i$ such that $e_i\Pi^0(Q) e_i$ has standard identity 
of degree $k$ (since some factor of $e_i\Pi^0(Q) e_i$ is algebra of matrices $\delta_i\times \delta_i$).
Moreover it implies that every polynomial map from $V$ to $\End_\Complex(V_i)$ which commutes with
all $G$-equivariant polynomial maps from $V$ to $\End_\Complex(V_i)$ accepts only scalar values thus
the inclusion of proposition \ref{prop0Cent} is in fact an isomorphism.
\begin{cor} \label{corZero}
Theorems \ref{thmCent} and \ref{thmPI} are valid when $\lambda=0$.
\end{cor}

\section{Regularity of the multiplication law}
Denote by $S_n$ the $\Complex^I$-bimodule $(\Sym^n(V^*)\otimes \End_\Complex(V_\Sigma))^G$, by $S$ the
graded algebra $(\Sym(V^*)\otimes \End_\Complex(V_\Sigma))^G$, by 
$T_n$ the $\Complex^I$-bimodule $(V^{*\otimes n}\otimes \End_\Complex(V_\Sigma))^G$ and by $T$ the graded
algebra $(T(V^*)\otimes \End_\Complex(V_\Sigma))^G$.
In this section we will show that all algebras of the family $\Pi^\lambda(Q)$ can 
be identified with an algebra which is $S$ as a vector space and the multiplication law in it
polynomially depends on $\lambda$.
For every $k=0,1,2,\dots$ we construct an operator 
\[
\pi^\lambda_k:T_k\To \bigoplus_{i=0}^k S_i
\]
such that
\begin{enumerate}
	\item $\pi^\lambda_k(x) = x$ for $x\in S_k$.
	\item $\pi^\lambda_k(x) \equiv x \mod{\delta\omega - \lambda}$ for any $x\in T_k$
	\item $\pi^\lambda_k(x_1 \omega x_2) = \pi^\lambda_{k-2}(x_1 \delta^{-1}\lambda x_2)$ for any $x_1\in T_i$ and $x_2\in T_j$ 
	with $i+j=k-2$.
	\item $\pi^\lambda_k(x)$ polynomially depends on $\lambda$.
\end{enumerate}
Then the family of operators $\pi^\lambda_k$ define an operator $\pi^\lambda$ acting from $T$ to $S$.
It is clear that $\pi^\lambda$ is a projection with image $S$, the second property of 
$\pi^\lambda_k$ guarantee that $\pi^\lambda(x)$ is equivalent to $x$ in algebra 
$\Pi^\lambda(Q)$, whilst the third property implies that equivalent in $\Pi^\lambda(Q)$
elements are mapped to identical elements. Combined this gives the isomorphism of $\Pi^\lambda(Q)$ and $S$
as filtered vector spaces and multiplication in $\Pi^\lambda(Q)$ transferred to $S$ can be easily written
as
\[
x\times y = \pi^\lambda(x\otimes y),
\]
which polynomially depends on $\lambda$. It is left to show that the family of operators with
properties (1)-(4) exist.

Clearly for $k=0$ and $k=1$ we may take an identity operators. Then we prove existence of 
$\pi^\lambda_k$ by induction. Fix some $\lambda\in \Complex^I$ and integer $k\ge2$. 
Define operators 
\[
\tau_i : T_k\oplus \bigoplus_{j=0}^{k-2} S_j \To T_k\oplus \bigoplus_{j=0}^{k-2} S_j \;\;\text{for $i=1,\dots, k-1$ as} 
\]
\begin{align*}
&\tau_i(x) = 0 \;\; \text{for}\; x\in \bigoplus_{j=0}^{k-2} S_j, \\
&\tau_i(x) = 0 \;\; \text{for}\; x\in T_k\; \text{such that $x$ is symmetric with respect to $i$-th and $i+1$-th arguments}, \\
&\tau_i(f \omega g) = f \omega g - \pi_\lambda^{k-2}(f\delta^{-1}\lambda g),
\end{align*}
which defines $\tau_i$ for $x\in T_k$ such that $x$ is antisymmetric with respect 
to $i$-th and $i+1$-th arguments. Put $\rho_i = 1-2 \tau_i$. We prove the following fact:
\begin{prop}
The family of operators $(\rho_i)$ satisfy conditions
\begin{enumerate}
\item $\rho_i^2=1$,
\item $\rho_i\rho_j = \rho_j\rho_i$ for $|i-j|>1$,
\item $\rho_i\rho_{i+1}\rho_i = \rho_{i+1}\rho_i\rho_{i+1}$,
\end{enumerate}
so $(\rho_i)$ induce a representation of the group of permutations of $k$ elements.
\end{prop}
\begin{proof}
Property (1) is easy. Consider the property (2). Assume $j>i$. It is enough to check the property 
for argument of the form
\[
x = f_1 \omega f_2 \omega f_3 \;\;\text{for}\; f_1\in T_{i-1},\; f_2\in T_{j-i-2},\; f_3\in T_{k-j-1}.
\]
Then
\[
\begin{split}
\rho_i \rho_j x - \rho_j \rho_i x = \pi_\lambda^{k-2}(f_1 \omega f_2 \delta^{-1}\lambda f_3 - 
f_1 \omega f_2 \delta^{-1}\lambda f_3) \\
= \pi_\lambda^{k-4}(f_1 \delta^{-1}\lambda f_2 \delta^{-1}\lambda f_3) - 
\pi_\lambda^{k-4}(f_1 \delta^{-1}\lambda f_2 \delta^{-1}\lambda f_3) = 0
\end{split}
\]
by the induction hypothesis. Consider the property (3). Denote by $\rho_i'$, $i=1,2,\dots,k-1$
the operator in $T_k$ which acts on $x\in T_k$ by interchanging of $i$-th and $i+1$-th arguments. Then,
clearly operators $\rho_i'$ satisfy conditions (1)-(3). Choose some $i\neq k-1$. Since there
is no element of $T_k$ which is antisymmetric with respect to arguments $i$, $i+1$, $i+2$ the
following operator vanishes on $T_k$:
\[
1-\rho_i'-\rho_{i+1}'-\rho_i'\rho_{i+1}'\rho_i'+\rho_i'\rho_{i+1}'+\rho_{i+1}'\rho_i' = 0.
\]
If we substitute $\rho_i' = 1-2 \tau_i'$ we obtain that 
\[
\tau_i' \tau_{i+1}' \tau_i' = \frac14 \tau_i' \;\;\text{and}\;\; 
	\tau_{i+1}' \tau_i' \tau_{i+1}' = \frac14 \tau_{i+1}'
\]
If we prove that
\[
\tau_i \tau_j \tau_i = \frac14 \tau_i \;\; \text{for} \; |i-j|=1,
\]
the property (3) would follow. But if $|i-j|=1$ then
\[
\tau_i \tau_j \tau_i = \tau_i^2 \tau_j \tau_i = \tau_i \tau_i' \tau_j' \tau_i' = \frac14 \tau_i \tau_i' = \frac14 \tau_i^2 = \frac14 \tau_i,
\]
here we consider $\tau_m'$, $m=1,2,\dots, k-1$ as an operator in $T_k\oplus \bigoplus_{i=0}^{k-2} S_i$ which
acts as zero on the component $\bigoplus_{i=0}^{k-2} S_i$ and use the 
equality $\tau_{m_1} \tau_{m_2} = \tau_{m_1} \tau_{m_2}'$ which is valid for 
$m_1, m_2 = 1,2,\dots,k-1$.
\end{proof}
Consider the representation of group of permutations of $k$ elements $\mathcal{S}_k$ 
given by operators $\rho_i$. 
Denote by $\bar{\varepsilon}$ the image of the element
\[
\varepsilon = \frac1{k!} \sum_{\sigma\in \mathcal{S}_k} \sigma
\]
of group algebra $\Complex \mathcal{S}_k$. Then we can expand every $\sigma$ as a product of
operators $\rho_i$, substitute $\rho_i = 1-2 \tau_i$ and represent
\begin{equation} \label{eqVarEps}
\bar{\varepsilon} = 1 + \sum_{i,j=1}^{k-1} \tau_i x_{ij} \tau_j \;\; \text{where $x_{ij}$ are some operators}.
\end{equation}
Then put $\pi_\lambda^k x = \bar{\varepsilon} x$ for $x\in T_k$. Check the required 
properties for $\pi_\lambda^k$. The property (1) follows from (\ref{eqVarEps}) and the fact
that all $\tau_i$ vanish on elements of $S_k$. Since all images of $\tau_i$ belong to the
ideal generated by $\delta\omega-\lambda$, the property (2) follows. The property (3) is true
since $\bar{\varepsilon} = \bar\varepsilon \rho_{i+1}$ implies $\bar\varepsilon = \bar\varepsilon (1 - \tau_{i+1})$ and 
\[
\bar\varepsilon (x_1 \omega x_2) = \bar\varepsilon(1 - \tau_{i+1})(x_1 \omega x_2) = \bar\varepsilon \pi_\lambda^{k-2} (x_1\delta^{-1}\lambda x_2) = 
\pi_\lambda^{k-2} (x_1\delta^{-1}\lambda x_2).
\]
The property (4) is obvious, so we have proved
\begin{prop}
The family of operators $\pi_\lambda^k$ satisfying properties (1) - (4) exist.
\end{prop}
The immediate corollary is
\begin{cor}\label{corS}
Every algebra $\Pi^\lambda(Q)$ is isomorphic as a filtered algebra to $S$ with multiplication
law $\Ast^\lambda$ which polynomially depends on $\lambda$ and is such that for 
any homogeneous $x$ of degree $i$ and homogeneous $y$ of degree $j$ the term
of degree $i+j$ in $x \Ast^\lambda y$ does not depend on $\lambda$.
\end{cor}

\section{Generic $\lambda$}
Due to corollary \ref{corS} we identify $\Pi^\lambda(Q)$ with $S$ with multiplication 
which depends on $\lambda$ polynomially. Denote this multiplication by $\Ast^\lambda$. 
Sometimes when $\lambda$ is fixed we will omit the sign $\Ast^\lambda$ and simply
write $x y$ instead of $x \Ast^\lambda y$ keeping in mind that the result depends on
$\lambda$ polynomially.
In this section we will prove the statement of theorems \ref{thmCent} and \ref{thmPI} for
some algebraically dense subset of the set 
\[
\eta = \{\lambda\in \Complex^I: \lambda \cdot \delta = 0\}.
\]
Namely, it will be proved for set where the following proposition holds:
\begin{prop} \label{propUnitExp}
There exist elements $f_1, \dots, f_n$ and $g_1, \dots, g_n$ in $S$ and rational functions
$\alpha_1, \dots, \alpha_n$ defined on $\eta$ such that
\[
\sum_{i=1}^n \alpha_i(\lambda) f_i \Ast^\lambda e_0 \Ast^\lambda g_i = 1
\]
for each $\lambda$ from some algebraically dense subset of $\eta$.
\end{prop}
\begin{proof}
It easily follows from the definition of deformed preprojective algebra that 
\[
\Pi^\lambda(Q)/\Pi^\lambda(Q) e_0 \Pi^\lambda(Q) \cong \Pi^{\lambda'}(Q'),
\]
where $Q'$ is the Dynkin diagram obtained from $Q$ by deleting vertex $0$ and $\lambda'$
is the restriction of $\lambda$ to vertices of $Q'$. It was proved in \cite{CrB} that 
deformed preprojective algebra of a Dynkin diagram is always finite dimensional and is zero
for all parameters except some number of hyperplanes. We will use the following implications:
\begin{enumerate}
\item the homogeneous subspace $S \Ast^0 e_0 \Ast^0 S$ of $S$ has finite codimension,
\item there exists $\lambda_0\in\eta$ such that $S \Ast^{\lambda_0} e_0 \Ast^{\lambda_0} S = S$.
\end{enumerate}
Choose some basis in $S \Ast^0 e_0 \Ast^0 S$ of the form $(a_i \Ast^0 e_0 \Ast^0 b_i)$ 
where $i$ ranges over the set of positive integers and all $a_i$ and $b_i$ are homogeneous elements
of $S$. It follows from the first statement that we can add some finite number of homogeneous 
elements of $S$ $x_1$, $x_2$, \dots, $x_n$ such that $x_i$ and $a_i \Ast^0 e_0 \Ast^0 b_i$
together form a basis of $S$. Now, for $\lambda \in \eta$ consider the set 
\[
B(\lambda) = \{x_i | i=1,\dots, n\} \cup \{a_i \Ast^\lambda e_0 \Ast^\lambda b_i | i=1, 2, \dots\}.
\]
It is again a basis of $S$ because each $a_i \Ast^\lambda e_0 \Ast^\lambda b_i$ equals to
sum of $a_i \Ast^0 e_0 \Ast^0 b_i$ and some terms of lower degree. Moreover every element of $S$
being expanded with respect to this basis has all coefficients polynomial in $\lambda$.

It follows from the statement (2) that there exist some $\lambda_0$ such that for $i=1,\dots, n$
\[
x_i = \sum_{k=1}^{K_i} f_i^k \Ast^{\lambda_0} e_0 \Ast^{\lambda_0} g_i^k
\]
where all $f_i^k$ and $g_i^k$ are elements of $S$.
Consider elements $y_i(\lambda)\in S$  for $i=1,\dots,n$ defined by 
\[
y_i(\lambda) = \sum_{k=1}^{K_i} f_i^k \Ast^\lambda e_0 \Ast^\lambda g_i^k.
\]
Consider an $n\times n$ matrix $Z(\lambda) = (z_{ij}(\lambda))$ where $z_{ij}(\lambda)$ is
the value of a coefficient near $x_i$ of the expansion of $y_j(\lambda)$ with respect to the basis
$B(\lambda)$. We have the following expansion of $y_j(\lambda)$ with respect to the basis $B(\lambda)$:
\[
\sum_{k=1}^{K_j} f_j^k \Ast^\lambda e_0 \Ast^\lambda g_j^k = \sum_{i=1}^n z_{ij}(\lambda) x_i +
\sum_{k=1}^{L_j} c_{jk}(\lambda) a_k \Ast^\lambda e_0 \Ast^\lambda b_k
\]
for some polynomial functions of $\lambda$ $c_{jk}(\lambda)$. Rewrite this as
\[
\sum_{i=1}^n z_{ij}(\lambda) x_i = \sum_{k=1}^{K_j} f_j^k \Ast^\lambda e_0 \Ast^\lambda g_j^k - 
\sum_{k=1}^{L_j} c_{jk}(\lambda) a_k \Ast^\lambda e_0 \Ast^\lambda b_k
\]
and consider it as a system of linear equations with indeterminates $x_1, \dots, x_n$. Clearly
it can be solved for such $\lambda$ that $\det{Z(\lambda)} \neq 0$ and the solution will depend on 
$\lambda$ rationally. If we expand $1$ with respect to the basis $B(\lambda)$ and then use 
this solution we obtain the required expansion. The set of $\lambda\in \eta$ for which 
$\det{Z(\lambda)} \neq 0$ is open. It is nonempty since $Z(\lambda_0)$ is the identity matrix, 
hence this set is dense. This completes the proof.
\end{proof}

Denote by $\eta'$ the subset of $\eta$ for which we proved the proposition above.
\begin{prop} \label{propCentGeneric}
For every $\lambda\in \eta'$ and every $x\in \mathcal{O}^\lambda(Q) = e_0 \Pi^\lambda(Q) e_0$ 
there exist $z(x)$ in the center of $\Pi^\lambda(Q)$ such that $e_0 z(x) e_0 = x$.
\end{prop}
\begin{proof}
Put 
\[
z(x) = \sum_{i=1}^n \alpha_i(\lambda) f_i x g_i.
\]
Then 
\[
e_0 z(x) e_0 = \sum_{i=1}^n \alpha_i(\lambda) e_0 f_i x g_i e_0 = 
\sum_{i=1}^n \alpha_i(\lambda) x f_i e_0 g_i e_0 = x
\]
since $\mathcal{O}^\lambda(Q)$ is commutative. Again, using commutativity of $\mathcal{O}^\lambda(Q)$
for any $y \in S$
\[
\begin{split}
y z(x) = \sum_{i=1}^n \alpha_i(\lambda) y f_i x g_i
= \sum_{i,j=1}^n \alpha_i(\lambda)\alpha_j(\lambda)f_j e_0 g_j y f_i x g_i
\\= \sum_{i,j=1}^n \alpha_i(\lambda)\alpha_j(\lambda)f_j x g_j y f_i e_0 g_i
= \sum_{j=1}^n \alpha_j(\lambda) f_j x g_j y = z(x) y.
\end{split}
\]
\end{proof}

\begin{prop}
For every $\lambda\in \eta'$ and every $q\in I$ the algebra $e_q \Pi^\lambda(Q) e_q$ has standard 
identity of degree $2 \delta_q$.
\end{prop}
\begin{proof}
For $x\in S$ construct a $n\times n$ matrix $M(x)$ over $\mathcal{O}^\lambda(Q)$ with elements
\[
m_{ij}(x) = \alpha_i(\lambda) e_0 g_i x f_j e_0.
\]
Then for $x, y \in S$ the matrix $M(x) M(y)$ has elements
\[
\begin{split}
\sum_{k=1}^n m_{ik}(x) m_{kj}(y) = 
\sum_{k=1}^n \alpha_i(\lambda) e_0 g_i x f_k e_0 
	\alpha_k(\lambda) e_0 g_k y f_j e_0 \\= 
	\alpha_i(\lambda) e_0 g_i x y f_j e_0 
	= m_{ij}(x y),
\end{split}	
\]
\[
\text{so}\;\; M(x y) = M(x) M(y).
\]
Denote by $p$ the matrix $M(1)$. Clearly $p$ is an idempotent and $M$ defines a homomorphism
from $\Pi^\lambda(Q)$ to $p \Mat(n, \mathcal{O}^\lambda(Q)) p$ where $\Mat(n, \mathcal{O}^\lambda(Q))$
denotes the algebra of $n \times n$ matrices over $\mathcal{O}^\lambda(Q)$. Construct an
inverse map $N: \Mat(n, \mathcal{O}^\lambda(Q)) \To S$. Let $A=(a_{ij})$ then put
\[
N(A) = \sum_{i,j=1}^n \alpha_j(\lambda) f_i a_{ij} g_j.
\]
Then we can check
\[
N(M(x)) = \sum_{i,j=1}^n \alpha_j(\lambda) f_i \alpha_i(\lambda) e_0 g_i x f_j e_0 g_j
=x \;\; \text{and}
\]
\[
m_{ij}(N(A)) = 
\sum_{k,l=1}^n \alpha_i(\lambda) e_0 g_i \alpha_l(\lambda) f_k a_{kl} g_l f_j e_0
\;\; \text{,which implies}
\]
\[
M(N(A)) = p A p.
\]
It proves that $M$ is an isomorphism. The algebra $\mathcal{O}^\lambda(Q)$ is a domain
(see \cite{CrB}). Hence it can be embedded into its field of fractions $F$. So the algebra
$p \Mat(n, \mathcal{O}^\lambda(Q)) p$ can be embedded into $p \Mat(n, F) p$ which is isomorphic
to $\Mat(r, F)$ where $r$ is the rank of $p$ in $\Mat(n, F)$. Denote by $p_q$ the matrix
$M(e_q)$ for $q\in I$. In a similar way $e_q \Pi^\lambda(Q) e_q$ can be embedded into 
$\Mat(r_q, F)$ where $r_q$ is the rank of $p_q$ in $\Mat(n, F)$. On the other hand 
$r_q = \tr{p_q}$ which is rational function of $\lambda$. Since $r_q$ can accept only a finite
number of values, namely $1, 2, \dots, n$ on the dense set $\eta'$ it is constant. 
In $\Pi^\lambda(Q)$ 
\[
\sum_{a\in \Arrows(Q)} [a, a^*] = \sum_{q\in I} \lambda_q e_q.
\]
Hence 
\[
\sum_{q\in I} \lambda_q r_q = \tr\sum_{q\in I} \lambda_q p_q = 0.
\]
Since this equality holds for all $\lambda$ from $\eta'$ which is dense in $\eta$
there is a constant $c\in \Complex$ such that $r_q = c \delta_q$ for $q\in I$. For $q=0$
\[
p_0 = M(e_0) = (\alpha_i(\lambda) e_0 g_i e_0 f_j e_0)
\]
so $p_0$ has rank $1$. It implies $c=1$ and $r_q = \delta_q$. We have proved that the algebra
$e_q \Pi^\lambda e_q$ for $\lambda\in \eta'$, $q\in I$ is isomorphic to some subalgebra of
the algebra of $\delta_q\times \delta_q$ matrices over the field $F$, so the standard identity
of degree $2 \delta_q$ is satisfied by Amitsur-Levitzki theorem.
\end{proof}

\section{Extending to the whole hyperplane}
To finish the proof of theorems \ref{thmCent} and \ref{thmPI} we need to make several steps.

\begin{prop} \label{propPreparePI}
For any $\lambda\in \Complex^I$ such that $\lambda \cdot \delta = 0$ and any $i\in I$
the algebra $e_i \Pi^\lambda(Q) e_i$ satisfies the standard identity of degree $2 \delta_i$.
\end{prop}
\begin{proof}
For $x_1, \dots x_{2\delta_i} \in e_i S e_i$ the sum
\[
\sum_{\sigma\in \mathcal{S}_{2 \delta_i}} \sign(\sigma) x_{\sigma(1)} \Ast^\lambda \dots 
\Ast^\lambda x_{\sigma(2 \delta_i)}
\]
is zero on an algebraically dense subset of $\lambda\in \Complex^I$, $\lambda \cdot \delta = 0$. Since 
it is polynomial in $\lambda$ it is zero for all $\lambda\in \Complex^I$, $\lambda \cdot \delta = 0$.
\end{proof}

\begin{prop} \label{propPrepareCent}
For every $\lambda\in \eta$ and every $x\in\mathcal{O}^\lambda(Q)$ there exist unique $z(x)$ in the
center of $\Pi^\lambda(Q)$ such that $e_0 z(x) e_0 = x$. 
\end{prop}
\begin{proof}
First note that if such a $z(x)$ exist then it is unique. Suppose the contrary. Then
there exists $a$ in the center of $\Pi^\lambda(Q)$ such that $e_0 a = 0$. Suppose $e_i a \neq 0$. Then
since $\Pi^\lambda(Q)$ is prime (see \cite{CrB}) there exist $y\in \Pi^\lambda(Q)$ such that
$e_0 y e_i a \neq 0$. Rewrite the last as $e_0 a y e_i$ and get a contradiction.

Then note that the degree of $z(x)$ is not greater then that of $x$. Let $z(x)'$ be the term of 
maximal degree of $z(x)$ and suppose that the degree of $z(x)'$ is greater then that of $x$. 
Clearly $z(x)'$ belongs to the center of $\Pi^0(Q)$, but $e_0 z(x)' e_0 = 0$ which contradicts
previous remark.

The algebra $\Pi^\lambda(Q)$ is finitely generated, and for any $x$ since the degree of $z(x)$
is bounded the problem of finding such $z(x)$ for any fixed $x$ is equivalent to some 
finite system of linear equations. Coefficients of the system depend on $\lambda$ polynomially.
Suppose the system has $m$ equations and $n$ indeterminates.
Consider the set $W$ of $\lambda$ for which the system has a unique solution. The system has
unique solution if and only if there exist equations $i_1, i_2, \dots, i_n$ in the system such
that the subsystem $i_1, i_2, \dots, i_n$ is nondegenerate (the set $U$ of $\lambda$ for which it is 
true is open) and the solution of equations $i_1, i_2, \dots, i_n$ satisfy other equations (the
set of $\lambda$ for which it is true is closed in $U$). Thus we obtain a sequence of open
sets $U_1, U_2, \dots, U_N$ and a sequence of sets $V_1, V_2, \dots, V_N$ each $V_i$ closed in
corresponding $U_i$. It follows that $W$ is covered by $U_1, U_2, \dots, U_N$ and intersection
of $W$ with each $U_i$ is closed. So $W$ is a closed set in the union of $U_1, U_2,\dots U_n$ hence
it is an intersection of some open set and some closed set.

Applying proposition \ref{propCentGeneric} together with the first remark in this proof we obtain 
that $W$ is an open set. Applying proposition \ref{prop0Cent} with first remark we obtain
that $W$ contains some neighbourhood of zero. So for any $x\in e_0 S e_0$ and any $\lambda$
there exist some constant $c\in \Complex$ such that there exist $z'(x)\in S$ which belongs to the center
of $\Pi^{c \lambda}(Q)$ and $e_0 z'(x) e_0 = x$. Let $x$ be a homogeneous element of degree $k$.
Define an operator $\phi$ on $T$ as a multiplication by
$c^{\frac{n}2}$ on each $T_n$. Then $\phi$ is an automorphism of algebra $T$ and maps $\delta \omega - c \lambda$
to $c \delta \omega - c \lambda$. It follows that $\phi(z'(x))$ belongs to the center of $\Pi^\lambda(Q)$
and $e_0 \phi(z'(x)) e_0 = c^{\frac{k}2} x$, so $z(x) = \phi(z'(x)) c^{-\frac{k}2}$ belongs to the center 
of $\Pi^{\lambda}(Q)$ and $e_0 z(x) e_0 = x$.
\end{proof}

\begin{proof}[Proof of the theorem \ref{thmCent}]
For any $\lambda\in \Complex^I$, $\lambda \cdot \delta=0$ take a map $\phi_\lambda$
from $\mathcal{O}^\lambda(Q)$ to the center of $\Pi^\lambda(Q)$ such that $e_0 \phi_\lambda(x) e_0 = x$ 
for all $x\in \mathcal{O}^\lambda(Q)$. By the proposition \ref{propPrepareCent} $\phi_\lambda$
is uniquely defined by this property so it is linear.
If $x, y\in \mathcal{O}^\lambda(Q)$ then $\phi_\lambda(x) \phi_\lambda(y)$ belongs to the
center of $\Pi^\lambda(Q)$ and $e_0 \phi_\lambda(x) \phi_\lambda(y) e_0 = x y$, so again by
the proposition \ref{propPrepareCent} $\phi_\lambda(x y) = \phi_\lambda(x) \phi_\lambda(y)$. 
Clearly $\phi_\lambda(e_0) = 1$. So $\phi_\lambda$ is a homomorphism.
The homomorphism $\phi_\lambda$ is an inclusion because for any $x\in \mathcal{O}^\lambda(Q)$
$x = e_0 \phi_\lambda(x) e_0$.

For any $i\in I$ put $\phi_\lambda^i(x) = e_i \phi_\lambda(x)$ for $x\in \mathcal{O}^\lambda(Q)$.
Then it is elementary to check that $\phi_\lambda^i$ is a homomorphism from algebra $\mathcal{O}^\lambda(Q)$ to
the center of $e_i \Pi^\lambda(Q) e_i$. It is an
inclusion because $\Pi^\lambda(Q)$ is prime (see \cite{CrB}), so if $x\neq 0$ belong to the center of 
$\Pi^\lambda(Q)$ then there exist $y\in\Pi^\lambda(Q)$ such that $e_i y x\neq 0$ hence $e_i x \neq 0$.

To prove that $\phi_\lambda^i$ is surjective suppose that $x$ belongs to the center of 
$e_i \Pi^\lambda(Q) e_i$, $x$ does not belong to the image of $\phi_\lambda^i$ and has 
the smallest possible degree. Let $x'$ be the term of highest degree of $x$ (we again identify
$\Pi^\lambda(Q)$ with $S$). Then $x'$ belongs to the center of $e_i \Pi^0(Q) e_i$ and thus 
there is homogeneous $y\in \mathcal{O}^\lambda(Q)$ such that $x' = \phi_0^i(y)$ (it already follows 
from the corollary \ref{corZero} that $\phi_0^i$ is surjective). 
Consider $z = \phi_\lambda(y)$ and $z'$~--- the term of the highest degree of $z$. Then $z'$ 
is in the center of $\Pi^0(Q)$ and $e_0 z' e_0$ is zero or equals to $y$. 
The first case is impossible due to the proposition \ref{propPrepareCent}. Thus 
$z' = \phi_0(y)$ and the term of maximal degree of $\phi_\lambda^i(y) = e_i z e_i$ 
equals to $x'$. It follows that $x - \phi_\lambda^i(y)$ has degree lower then $x$ and does
not belong to the image of $\phi_\lambda^i$ thus obtaining a contradiction.
\end{proof}

\begin{proof}[Proof of the theorem \ref{thmPI}]
The statement of the theorem \ref{thmPI} follows from the proposition \ref{propPreparePI}
and the fact that if $k$ is such that for any $x_1, x_2, \dots, x_{k} \in e_i S e_i$
\[
\sum_{\sigma\in \mathcal{S}_{k}} \sign(\sigma) x_{\sigma(1)} \Ast^\lambda \dots 
\Ast^\lambda x_{\sigma(k)} = 0,
\]
then denoting by $x_i'$ the term of maximal degree of $x_i$ we get
\[
\sum_{\sigma\in \mathcal{S}_{k}} \sign(\sigma) x'_{\sigma(1)} \Ast^0 \dots 
\Ast^0 x'_{\sigma(k)} = 0,
\]
so from the corollary \ref{corZero} $k\ge 2 \delta_i$.
\end{proof}

\section{Proof of the theorem \ref{thm1}}

Consider a quiver $C_n$ with $n$ vertices $I=\{1, 2, \dots, n\}$ which form a chain:
\[
\begin{tiny}
	\xymatrix{
	n & n-1 \ar@{->}[l]^{a_{n-1}} & n-2 \ar@{->}[l]^{a_{n-2}} & \cdots \ar@{->}[l] & 1 \ar@{->}[l]^{a_1}
	}
\end{tiny}	
\]
Suppose we have a sequence of complex numbers $\lambda = (\lambda_i)$, $i=1, \dots, n-1$. Consider an 
algebra 
\[
R^\lambda_n = e_n \left( \Complex \bar{C_n} / (\sum_{i=1}^{n-2} [a_i, a_i^*] - a_{n-1}^*a_{n-1} - \sum_{i=1}^{n-1} \lambda_i e_i) \right) e_n.
\]
\begin{prop}\label{propPreparethm1}
The algebra $R^\lambda_n$ is isomorphic to the algebra $\Complex[x]/P(x)$ via an isomorphism
sending $x$ to $a_{n-1} a_{n-1}^*$ where P(x) is
a polynomial given by 
\[
P(x) = x(x+\lambda_{n-1})(x+\lambda_{n-1}+\lambda_{n-2})\dots(x+\sum_{i=1}^{n-1} \lambda_i).
\]
\end{prop}
\begin{proof}
If $n=1$ both algebras are isomorphic to $\Complex$. We proceed by induction. For $n>1$ the
algebra $R^\lambda_n$ splits as a vector space:
\[
R^\lambda_n = \Complex \oplus a_{n-1} e_{n-1} \left( \Complex \bar{Q} / (\sum_{i=1}^{n-1} [a_i, a_i^*] - a_{n-1}^*a_{n-1} - \sum_{i=1}^{n-1} \lambda_i e_i) \right) e_{n-1} a_{n-1}^*.
\]
Then,
\[
\begin{split}
e_{n-1} \left( \Complex \bar{C_n} / (\sum_{i=1}^{n-1} [a_i, a_i^*] - a_{n-1}^*a_{n-1} - \sum_{i=1}^{n-1} \lambda_i e_i) \right) e_{n-1}
 \\\cong (R^\lambda_{n-1} \ast \Complex[a_{n-1}^* a_{n-1}]) / (a_{n-2} a_{n-2}^* - a_{n-1}^*a_{n-1} - \lambda_{n-1}e_{n-1}),
\end{split}
\] 
where we denote by $\ast$ the free product of algebras. By the induction hypothesis the last is isomorphic to
\[
(\Complex[a_{n-2} a_{n-2}^*]/P^-(a_{n-2} a_{n-2}^*) \ast \Complex[a_{n-1}^* a_{n-1}]) / (a_{n-2} a_{n-2}^* - a_{n-1}^*a_{n-1} - \lambda_{n-1}e_{n-1})
\]
for 
\[
P^-(x) = x(x+\lambda_{n-2})(x+\lambda_{n-2}+\lambda_{n-3})\dots(x+\sum_{i=1}^{n-2} \lambda_i),
\]
so
\[
\begin{split}
e_{n-1} \left( \Complex \bar{C_n} / (\sum_{i=1}^{n-1} [a_i, a_i^*] - a_{n-1}^*a_{n-1} - \sum_{i=1}^{n-1} \lambda_i e_i) \right) e_{n-1}
 \\\cong \Complex[a_{n-1}^* a_{n-1}]/P^-(a_{n-1}^* a_{n-1} + \lambda_{n-1}),
\end{split}
\] 
therefore
\[
R^\lambda_n \cong \Complex[a_{n-1} a_{n-1}^*]/(P^-(a_{n-1} a_{n-1}^* + \lambda_{n-1})a_{n-1} a_{n-1}^*)
\]
and it can be easily seen that
\[
P^-(a_{n-1} a_{n-1}^* + \lambda_{n-1})a_{n-1} a_{n-1}^* = P(a_{n-1} a_{n-1}^*).
\]
\end{proof}

The theorem is now valid because $e_c\Pi^\lambda(Q)e_c$ defined as in the statement of the theorem
is isomorphic to the free product
of algebras $R_{\deg P_i - 1}^{\lambda^i}$ factored by relation
\[
\sum_{i=1}^n a_{i1} a_{i1}^* = \mu e_c,
\]
where 
\[
\lambda^i = (\alpha_{i \deg P_i - 2} - \alpha_{i \deg P_i - 1}, \dots, \alpha_{i 1} - \alpha_{i 2}, -\alpha_{i 1})
\]
and by the proposition \ref{propPreparethm1} each $R_{\deg P_i - 1}^{\lambda^i}$ is isomorphic 
to 
\[
\Complex[a_{i1} a_{i1}^*]/P_i(a_{i1} a_{i1}^*).
\]


\begin{thebibliography}{x}
\bibitem{ours} M. Vlasenko, A. Mellit, Yu. Samoilenko. On algebras generated by
linearly related generators with given spectrum, to appear.
\bibitem{CrB} W. Crawley-Boevey, M. P. Holland. Noncommutative deformations of Kleinian
singularities, Duke Math. J., 92 (1998), 605-635.
\bibitem{Rowen1990} Rowen L.H. Ring theory. Academic Press, 1991.
\bibitem{McKay} J. McKay. Graphs, singularities, and finite groups, Proc. Sympos.
Pure Math., 37 (1980), 183-186.
\end{thebibliography}
\end{document}